\documentclass{article}
\usepackage{amssymb}

\newtheorem{theorem}{Theorem}

\newtheorem{conjecture}[theorem]{Conjecture}

\newtheorem{definition}[theorem]{Definition}

\newenvironment{proof}[1][Proof]{\noindent\textbf{#1.} }{\ \rule{0.5em}{0.5em}}
\input{tcilatex}

\begin{document}

\title{Spiral Chains: The Proofs of Tait's and Tutte's Three-Edge-Coloring
Conjectures }
\author{I. Cahit}
\date{}
\maketitle

\begin{abstract}
In this paper we have shown without assuming the four color theorem of
planar graphs that every (bridgeless) cubic planar graph has a
three-edge-coloring. This is an old-conjecture due to Tait in the squeal of
\ efforts in settling the four-color conjecture at the end of the 19th
century. We have also shown the applicability of our method to another
well-known three edge-coloring conjecture on cubic graphs. Namely Tutte's
conjecture that "every 2-connected cubic graph with no Petersen minor is
3-edge colorable". Hence the conclusion of this paper implies another
non-computer proof of the four color theorem by using spiral-chains in
different context.
\end{abstract}

\bigskip

\begin{center}
\bigskip \textit{To the memory of W.T. Tutte}
\end{center}

\section{Introduction}

Assuming the four color theorem (4CT) it is easy to show that every cubic
(bridgeless) planar graph has a three-edge-coloring since Tait found an
equivalent formulation of the 4CT in terms of 3-edge-coloring.\footnote{%
Assume that the regions of a cubic planar graphs colored properly by red
(R), blue (B), green (G) and yellow (Y) colors. Then if two neighbor regions
are colored by RG or GY then color its boarder (edge in the cubic graph) by
G, if two neighbor regions are colored by RG or BY then color it boarded by
R and if two neighbor regions are colored by RY or BG then color its boarder
by B. Then the theorem follows.}. This is an old-conjecture due to Tait in
the squeal of \ efforts in settling the four-color conjecture at the end of
the 19th century. In this paper we have given an independent proof of this
famous conjecture. Again as we have recently introduced the new concept
"spiral chains" in maximal planar graphs in the new non-computer proof of
the four color theorem, we continue the use of spiral chains in the cubic
plane graph for its three-edge-coloring [15]. Note that there are other
equivalent formulation of 4CT in terms of three edge-coloring of the pairs
of binary trees [4].[5]. We have first shown that any cubic planar graph can
be decomposed into vertex disjoint spiral chains which in turn expressed as
the union of \textit{folded\footnote{%
Some of the edges are called "hairs" used twice, where \textit{hair}-edge is
the edge of the comb-caterpillar other than the backbone-edge (an edge
having two ends degree three).}} special spiral caterpillar (spiral-combs).
Then based on the spiral-comb(s) decomposition we have provided a simple
three-edge-coloring algorithm for any given cubic planar graph. In the
second part of the paper we will show that the same method can also be
applied to a similar but more general three-edge cubic graph (non-planar)
coloring conjecture of Tutte which has been recently settled by a series
lengthy papers using the classical four color theorem proof technique. That
is

\begin{conjecture}
\bigskip Every $2$-connected cubic graph with no Petersen minor is 3-edge
colorable.
\end{conjecture}

This extends the four color theorem as Tait showed that the four color
theorem is equivalent to the statement that every planar $2$-connected cubic
graph is $3$-edge colorable. It also implies that certain non-planar graphs
are 3-edge colorable. Let us say that G is \textit{apex }if $G\diagdown v$
is planar for some $v$ and $G$ is \textit{doublecross} if it can be drawn in
the plane with crossings, but with atmost two crossings and with all the
crossings on the boundary of the finite region. Both apex and doublecross
graphs have no Petersen minor so Tutte's conjecture implies:

\begin{conjecture}
(RST) (1) Every 2-connected apex cubic graph is 3-edge colorable and (2)
every 2-connected doublecross cubic graph is 3-edge colorable.
\end{conjecture}

\bigskip The above conjecture is the outcome of joint efforts of Robertson,
Seymour and Thomas (RST) [6,7,8] and it appears that excluding Roberstson,
Sanders joins in the team for edge 3-coloring of cubic apex and doublecross
graphs [9,10].

Tutte and RST equivalences has been shown in [6] that "a minimum
counter-example" mean a 2-connected cubic graph $G$ with no Petersen minor
which is not $3$-edge colorable, with $\ |V(G)|$ minimum. That is every
minimum counter-example is either apex or doublecross. Eventually Tutte's
conjecture has been settled by a series of papers [6]-[10] by using the
classical proof method of the four color theorem [12],[13],[14]. The concept
of spiral chain has been introduced by the author for an non-computer proof
of the four-color theorem in 2004 [15]. We have applied spiral-chain
coloring technique to the maximal planar graphs and we will use it the same
for cubic planar and almost planar (apex and doublecross) cubic graphs.
Therefore the main aim of this paper is to give an independent proof of the
famous Tait's conjecture for bridgeless cubic planar graphs by using
spiral-chain edge-coloring and then to show that it would be possible to
extend it to the proof of conjecture RST. Hence to the proof of Conjecture 1.
\begin{verbatim}
 
\end{verbatim}

\section{Spiral-chains in cubic graphs}

Throughout the paper the numbers $1,2,3$ attached on the edges of a graph
correspond, respectively to the colors (R)ed, (O)range and (G)reen\footnote{%
These selection of colors may be further informational in the description of
spiral-chain edge coloring in the cubic graph when they are viewed as the
colors of the "traffic lights". That is Red for STOP, Orange for READY and
Green for GO.}. That is $C=\{R,O,G\}\equiv \{1,2,3\}$. Let $G=G_{0}$ be a $2$%
-connected (bridgeless) cubic planar graph drawn in the plane without
crossings. Let $B$ be denote its outer-boundary cycle consist of vertices $%
v_{1},v_{2},...,v_{k},k\geq 3$. Without loss of generality select the edges
starting from vertex $v_{1}$ends at vertex $v_{k}$ in the
clockwise-direction (or anti-clockwise-direction) in the plane. That is we
have selected the path from vertex $v_{1}$ to $v_{k}$ on the $B$. We can do
this always by selecting unused right-most edge in the plane, that is if the
vertex $v_{i}$ is in the spiral $S$ then $v_{i+1}$ will be in the $S$ only
if $(v_{i},v_{i+1})$ is the right-most edge incident at $v_{i}$. Delete
these selected vertices $v_{1},v_{2},...v_{k}$ from the boundary cycle of $%
G_{0}$ and obtain the subgraph $G_{1}=G_{0}\diagdown
\{v_{1},v_{2},...,v_{k}\}.$ Note that we have not selected the edge $%
(v_{k},v_{1})$ since $(v_{k},v_{1})\cup P(v_{1},v_{k})$ is the
outer-boundary cycle of $G_{0}$. Then continue in this way for the other
(inner) vertices on the outer-boundary vertices of the subgraph $G_{1},$
that is select the edge $(v_{k},v_{k+1}),$ where $\ v_{k}\in V(G_{0})$ and $%
v_{k+1}\in V(G_{1})$ and the outer-boundary vertices $%
v_{k+1},v_{k+2},...,v_{k+p}$ of $G_{1}$ in the clockwise direction and
delete them form $G_{1}$ and the repeat the same process for $%
G_{2}=G_{1}\diagdown \{v_{k+1},v_{k+2},...,v_{k+p}\}$. Let us assume that
repeatation of this spirally selection of the vertices results an empty
graph at stage $s$. Now union of deleted vertices of $%
G_{0},G_{1},G_{2},...,G_{s-1}$ forum the spiral-chain $S$ of $G.$We may
describe spiral-chain $S$ in terms of consecutive edges between the deleted
edges as $\{e_{1},e_{2},e_{3},...,e_{k-1}\}\cup (e_{k})\cup
\{e_{k+1},e_{k+2},...,e_{k+p-1}\}\cup ...,$where $e_{k}=$ $(v_{k-1},v_{k+1})$
is the edge connecting two spiral segments of $G_{0}$ and $G_{1}$. Not all $%
2 $-connected cubic planar graphs are hamiltonian, see for example famous
Tutte graph in Fig.2. Otherwise the proof of Tait's conjecture would be
immediately obtained by coloring hamiltonian cycle with colors 1 and 2 and
coloring other edges with color 3. But deciding whether a cubic planar graph
has an hamiltonian cycle or not is not an easy task. Therefore on Tait
coloring of cubic planar graphs hamiltonicity is not our primary concern. On
the other hand a $2$-connected cubic planar graph may have a hamiltonian
cycle but have several spiral chains (see Fig. 1 and 5). Edge coloring of
these figures will be explained later.

Let $T$ be tree with $n$ vertices. Base tree $T_{b}$ of $T$ is obtained by
removing end-vertices of $T$. A tree $T$ is called \textit{caterpilla}r if
its base tree is a path. We will use in our spiral-chain edge-coloring a
special caterpillar called comb-tree:

\begin{definition}
\bigskip Comb tree $T_{c}$ is a caterpillar in which all vertex degrees are $%
3$ or $1.$
\end{definition}

\begin{center}
\bigskip \FRAME{ftbpFU}{331.5625pt}{331.4375pt}{0pt}{\Qcb{(a) Coloring by
the use of an hamiltonian cycle, (b) spiral-chain edge-coloring}}{}{Figure}{%
\special{language "Scientific Word";type "GRAPHIC";maintain-aspect-ratio
TRUE;display "USEDEF";valid_file "T";width 331.5625pt;height
331.4375pt;depth 0pt;original-width 569.5625pt;original-height
569.375pt;cropleft "0";croptop "1";cropright "1";cropbottom "0";tempfilename
'IJ4GNL00.wmf';tempfile-properties "XPR";}}
\end{center}

We will be investigating 3-edge coloring of $2$-connected cubic planar
graphs under two classes of cubic planar graphs since the "triangles" (cycle
of length three) are the main issue in the uniqueness of the three colorings
in planar graphs. For example T. Fowler and R. Thomas have settled a
conjecture of Fiorini and Wilson, and Fisk about uniquely colorable planar
graphs e.g.,

\begin{theorem}
A 3-regular planar graph has a unique edge 3-coloring if and only it can be
obtained from the complete graph on four vertices by repeatedly replacing
vertices by triangles.
\end{theorem}

We will return to the uniquely edge 3-coloring planar graphs from the
point-view of spiral-chain coloring.

\section{Spiral-chain edge three-coloring}

\subsection{Triangle-free cubic bridgeless planar graphs}

Let us start with a result on three edge-coloring of cubic planar graphs
with triangles.

\begin{theorem}
Three edge-coloring of a cubic planar graph $G_{\vartriangle }$ with
triangles can be obtained from three edge-coloring of a triangle-free cubic
bridgeless planar graph $G$ .
\end{theorem}

\begin{proof}
Let $G_{\vartriangle }$ be denote a cubic (bridgeless) planar graph. Then $%
G_{\vartriangle }$ can only contains (a) vertex disjoint triangles, (b) two
triangles with only one edge common, call it \textit{twin-triangle} and (c)
if it contains three mutually adjacent triangles then $G_{\vartriangle
}=K_{4}.$Assume that all triangle-free cubic planar graphs have been three
edge-colored. Consider any cubic graph $G_{\vartriangle }$ which contains an
triangle $C_{3}$ without an common edge with another triangle. That is let $%
V(C_{3})=\{v_{1},v_{2},v_{3}\}$ and $E(C_{3})=%
\{(v_{1},v_{2}),(v_{2},v_{3}),(v_{3},v_{1})\}$. $G_{\vartriangle }$ must
have also vertices $v_{4},v_{5},v_{6}$ such that $\left( v_{1},v_{4}\right)
,\left( v_{2},v_{5}\right) ,\left( v_{3},v_{6}\right) \in E(G_{\vartriangle
})$ where $v_{4}\neq v_{5}\neq v_{6}$. Let $C=\{i,j,k\},i\neq j\neq k$ be
the set of the three colors, i.e., $i,j,k\in \{1,2,3\}$. In any proper
edge-coloring of $G_{\vartriangle }$ the subgraph formed by the triangle $%
C_{3}$ together with the edges $\left( v_{1},v_{4}\right) ,\left(
v_{2},v_{5}\right) ,\left( v_{3},v_{6}\right) $ must have an unique
edge-coloring as shown in Fig. 2(a). Since the colors of the edges $\left(
v_{1},v_{4}\right) ,\left( v_{2},v_{5}\right) ,\left( v_{3},v_{6}\right) $
are $i,j,k$ we can shrink the triangle $C_{3}$ into a vertex $u_{1}$ and
obtain a new cubic graph (see also Fig. 2(a)). We repeat this process for
every triangle till we arrive to a cubic-graph without distinct triangles.
Order the new vertices so created as $u_{1},u_{2},...,u_{k}.$ If $%
G_{\vartriangle }$contains twin-triangle then it must have vertices $%
v_{1},v_{2},v_{3},v_{4}$ and edges $%
(v_{1},v_{2}),(v_{2},v_{3}),(v_{2},v_{4}),(v_{1},v_{3}),(v_{3},v_{4})$ and
two more edges $(v_{1},v_{5})$ and $(v_{4},v_{6})$ to make vertex degrees of
the twin-triangle all three. Now all possible three edge-coloring of the
subgraph has been shown in Fig. 2(b). Note that we have $v_{5}\neq v_{6}$
otherwise it implies a bridge in $G_{\vartriangle }$. Since the color of the
edges $(v_{1},v_{5})$ and $(v_{4},v_{6})$ is $i\in \{1,2,3\}$ we delete the
twin-triangle from $G_{\vartriangle }$ and join $v_{5}$ and $v_{6}$ with an
new edge $e_{1}=(v_{5},v_{6})\,$to obtain another possibly triangle free
cubic planar graph. If not we repeat the above processes for each
twin-triangles in $G_{\vartriangle }$ till we get an triangle-free cubic
planar graph $G$. Let $e_{1},e_{2},...,e_{q}$ be set of new edges inserted
during the deletion of twin-triangles. Since we have assumed that $G$ is a
three edge-colorable we start from a coloring of $G$ and re-insert all the
deleted triangles for each $u_{i}\in V(G),i=1,2,...,k$ and twin-triangles
for each $e_{i}\in E(G),i=1,2,...,q$ back into $G$ but this time with
suitable edge-colors based on the coloring of $G$ and the cases in Fig. 2.
Hence we have obtained three edge-coloring of $G_{\vartriangle }$. Note that 
$K_{4}$ is an isolate case and is shown in Fig. 2(c).
\end{proof}

\begin{center}
\FRAME{ftbpFU}{269.9375pt}{190.5pt}{0pt}{\Qcb{Deletion of triangles from
cubic planar graphs.}}{}{Figure}{\special{language "Scientific Word";type
"GRAPHIC";maintain-aspect-ratio TRUE;display "USEDEF";valid_file "T";width
269.9375pt;height 190.5pt;depth 0pt;original-width 892.125pt;original-height
628pt;cropleft "0";croptop "1";cropright "1";cropbottom "0";tempfilename
'IIRF8M02.wmf';tempfile-properties "XPR";}}
\end{center}

\subsection{Uniquely 3-edge colorable planar graphs}

Characterization of uniquely 3-edge colorable graphs have been given by
Towler and Thomas. Since these cubic planar graphs can only be obtained from
the complete graph on four vertices by repeatedly replacing vertices by
triangles, we note that they are hamiltonian and admit a very simple unique $%
3-$edge coloring e.g., consider the coloring of Fig. 2(a). Let $G_{u}$ be an
uniquely 3-edge colorable graph. Let the vertices $(v_{i},v_{j},v_{k})$
induce an triangle in $G_{u}$ and let $C=\{1,2,3\}$. We say the edge colors
of \ an triangle is positive "+" if the colors of the edges $%
(v_{i},v_{j}),(v_{j},v_{k}),(v_{k},v_{i})$ respectively are $\ 1,2,3$ and
negative "-" if the colors of its edges respectively are $2,1,3$. That is in
a positively signed triangle, numbers corresponding colors are rotated in
clockwise direction while negatively signed triangle numbers are rotated in
anti-clockwise direction (see Fig. 5). We also say that a vertex $v$ in $%
G_{u}$ is mature if it is on the corner of an triangle otherwise we call
immature vertex. If all vertices of $G_{u}$ are mature we say the graph is
complete and denote it by $K_{u},$ i.e., complete with respect to the number
of triangles$.$ For example the cubic graph shown in Fig. 5 is a complete
graph on six triangles. Some properties of $G_{u}$ are:

- Except $K_{4}$ (complete graph with four vertices), for any cubic graph $%
G_{u}$ no two triangles have a common edge.

- \ \ Let $S$ be a spiral-chain in $G_{u}.$ For an 3-edge coloring of $G_{u}$
the sign of all triangles incident to outer-boundary cycle have the same
polarity, say all are +. Furthermore because of the character of the
construction of uniquely 3-edge colorable graphs they constitute
shell-structure (known also anti-matroids) of triangles in the graph and
triangles in the consecutive shells must be labeled with opposite signs. For
example in Fig. 5 all triangles labeled with positive signs are the
triangles in the outer shell and neighbor shell triangles are labeled with
negative signs. Since spiral-chains are shelling structure [15], we can use
this property when we color the edges of $G_{u}$. That is we can extract the
sign of the triangle from the knowledge of the location of the edge in the
spiral chain and from the sign of the previous triangle when we meet one
along the edges of the spiral chain. We omit details of the spiral-chain
coloring algorithm for 3-edge uniquely colorable graphs since the algorithm
given in the next section may as well be applied for this class of cubic
planar graphs.

\subsection{Spiral-chain edge-coloring}

In this section we will give an three edge-coloring algorithm which uses
spiral chains of the given cubic bridgeless planar graph. Let \thinspace $%
C=\{R,O,G\}.$ In the figures we use numbers$1,2,3$ on the edges instead of
letters to avoid confusion. Let $S$ be denote a spiral-chain of $G$. As we
go through the vertices of $S$ we construct a path $P_{s}$ with vertices $%
v_{1},v_{2},...,v_{k}$ and assume that $\mid V(P)\mid =\mid V(G)\mid ,$i.e., 
$k=n.$ We will investigate cubic graphs with several spiral-chains later.
Edges incident to $P_{s}$ together with the edges $E(P_{s})$ of the path
constitute a comb-caterpillar tree $C_{s}(V_{s},E_{s}),$ where $%
V_{s}=\{v_{1},v_{2},...,v_{k}\}$ and $E_{s}=E(P_{s})\cup \{\overline{E(P_{s}}%
):$ $\{v_{i},v_{j}\}\mid $ $v_{i},v_{j}\in V(P_{s})$ and $d(v_{i},v_{j})\geq
2\}$. An edge $b$ of the $C_{s}$ is called the \textit{backbone}-edge if $\
b\in E(P_{s})$ and is called the \textit{hair}-edge $h$ if $h\in \overline{%
E(P_{s}})$ . It is clear that $C_{s}=G$ but some of the edges of $\overline{%
E(P_{s}})$ are counted twice since any hair-edge of $C_{s}$ has two ends in
the spiral-chain $S.$ Again as we go through the edges of $S$ we call that
hair edge $h\in \overline{E(P_{s}})$ is \textit{real }if it is encountered
for the first-time and \textit{imaginary} if it is encountered for the
second-time. In the figures all the real-hair edges have been shown in
straight-lines while the imaginary-hair edges are shown with dashed-lines.
The importance of the real and imaginary hairs will be clearer when deciding
the color of the current edge in the spiral chain. We have classify the
colors into two classes. The color green $G$ is called as the "primary"
color while the colors orange $O$ and red $R$ are called the "secondary"
colors. In the coloring algorithm given below primary color $G$ is used
mostly for the spiral-chain backbone edges and the secondary colors $O$ and $%
\ R$ are used mostly for the hair-edges.

\bigskip

\textbf{Algorithm A. }\ \textit{Spiral-chain edge-coloring (for single
spiral-chain)}

Let the cubic planar graph G be drawn in plane without any crossings. As the
first step construct in clockwise direction a spiral-chain $S$ of $G$
starting from a vertex $v_{1}$ of its outer-boundary cycle. Let $C_{s}$ be
the corresponding comb of $S$. Start coloring the edges of the backbone
edges of $C_{s}$ alternatingly by green and orange colors i.e., $(G,O)$%
-chain and the hair edges of $C_{s}$ by the red $R$ color. \ If all vertices
of $S$ is covered and edges incident to vertices are properly colored by $%
G,O,$ and $R$ then the coloring is the desired proper three edge-coloring.

a) If for some backbone edge $e_{i}=(v_{i},v_{i-1})$ is forced to be colored
by $R$ because of the color (it must be orange $O$ since $G$ is primary
color) of the hair $h_{i-1}$ at $v_{i-1}.$ Furthermore if the hair $h_{i}$
had been colored by $R$ before then both $e_{i}$ and $h_{i}$ would be $R$.
In this case there must be a vertex $v_{j},$ with $j>i$ and a spiral
sub-path $P(v_{r},v_{s}),$ with $s>r<j>i$. But the path $%
P(v_{i},v_{j})=(v_{i},v_{r})\cup P(v_{r},v_{s})\cup (v_{s},v_{j})$ is a $%
(R,O)-$Kempe chain, where $(v_{i},v_{r})$ and $(v_{s},v_{j})$ are the
hair-edges. Therefore we perform Kempe-chain color switching for the edges
of $P(v_{i},v_{j})$ and resolves the problem at vertex $v_{i}$.

a) If for some backbone edge $e_{i}=(v_{i},v_{i-1})$ is forced to be colored
by $O$ because of the color of the hair $h_{i-1}$ at $v_{i-1}.$ Furthermore
if the hair $h_{i}$ had been colored by $O$ before then both $e_{i}$ and $%
h_{i}$ would be $O$. In this case there must be a vertex $v_{j},$ with $j>i$
and a spiral sub-path $P(v_{r},v_{s}),$ with $s>r<j>i$. But the path $%
P(v_{i},v_{j})=(v_{i},v_{r})\cup P(v_{r},v_{s})\cup (v_{s},v_{j})$ is a $%
(O,G)-$Kempe chain, where $(v_{i},v_{r})$ and $(v_{s},v_{j})$ are the
hair-edges. Therefore we perform Kempe-chain color switching for the edges
of $P(v_{i},v_{j})$ and resolves the problem at vertex $v_{i}$.

That is in critical cases we have always suitable Kempe-chain with both ends
are hairs of the spiral comb-tree. Hence we can perform Kempe-chain color
switching within the spiral-chain of $G$ when we needed to resolve color
duplication at vertex $v_{i}.$

Fig. 1(b) and Fig. 3(b) illustrate the spiral-chain edge-coloring algorithm.

Algorithm A can be easily modified to cubic planar graphs having more than
one spiral chains. In this case if the current vertex, say $v_{a}$ of
spiral-chain $S_{i}$ adjacent to two vertices $\ v_{b},v_{c}$ that have
already been contained in $S_{i}$ before or in some other spiral-chain $%
S_{j},j<i$ then we cannot go further. Hence we select closest new vertex
incident $S_{k},k=1,2,..,i$ and start new spiral-chain $S_{i+1}$ from that
vertex. When all vertices of $G$ covered we obtain vertex disjoint set of
spiral-chains $S_{1},S_{2},...,S_{m}.$ Note that some of spiral-chains
degenerate to a single vertex. One can easily notice that spiral-chains $%
S_{1},S_{2},...,S_{m}$ are topologically nested in the plane. See for
example spiral-chains in Fig. 5 and 6. Now we ready to give spiral-chain
edge-coloring algorithm for this case:

\bigskip

\textbf{Algorithm B. }Spiral-chain edge-coloring (for multiple
spiral-chains).

Let $S_{1},S_{2},...,S_{m}$ be the set of spiral-chains of the cubic planar
graph. Apply Algorithm A to each spiral-chain $S_{i},i=1,2,...,m$. Note that
for each application of Algorithm A we select always the color "green" $\ G$
as the primary color. This selection would minimize execution of number of
the steps (a) and (b) in Algorithm A. That is it will color in green,
maximum number of the backbone edges in the spiral-chains.

\FRAME{ftbpFU}{191.625pt}{331.4375pt}{0pt}{\Qcb{(a) Spiral-chain and
spiral-cut, (b) and its spiral three edge-coloring.}}{}{Figure}{\special%
{language "Scientific Word";type "GRAPHIC";maintain-aspect-ratio
TRUE;display "USEDEF";valid_file "T";width 191.625pt;height 331.4375pt;depth
0pt;original-width 486.0625pt;original-height 844.1875pt;cropleft
"0";croptop "1";cropright "1";cropbottom "0";tempfilename
'IJ6C1K01.wmf';tempfile-properties "XPR";}}\FRAME{ftbpFU}{141.25pt}{%
382.4375pt}{0pt}{\Qcb{Spiral comb tree in the coloring algorithm.}}{}{Figure%
}{\special{language "Scientific Word";type "GRAPHIC";maintain-aspect-ratio
TRUE;display "USEDEF";valid_file "T";width 141.25pt;height 382.4375pt;depth
0pt;original-width 319.0625pt;original-height 872.5625pt;cropleft
"0";croptop "1";cropright "1";cropbottom "0";tempfilename
'IJ1OWP04.wmf';tempfile-properties "XPR";}}

\section{Tutte's three edge-coloring conjecture}

The main idea of this section is to show that possible cubic graph minimal
counter-examples in Conjecture 2 (RST) can be three edge-colorable by using
spiral-chain edge coloring algorithm that has been given in the previous
section. In other words the contribution of this section is a finishing
touch to the deep result of Robertson, Seymour and Thomas [6] that the only
possible minimal counter-examples to Tutte's conjecture (Conjecture 1) would
be apex and doublecross cubic non-planar graphs. We only need to give a
suitable modification in the construction of spiral chains since the
possible counter-examples are specific non-planar cubic graphs i.e., apex
and doublecross. Let $G_{ad}$ denote apex or doublecross cubic graph.

Let us note the following observation:

In the previous section we have agreed that spiral-chains in the cubic
planar graphs are rotated in the clockwise or anti-clockwise. One of the
reason of this is to \textit{not} loss our way in the edge-coloring
algorithm in the cubic planar graph and ordering the spiral-chains in an
nested fashion. In any cubic non-planar graph $G_{ad}$ planarity are
destroyed either by a special vertex or two crossings of a pair of edges
within the same face. Suppose we start construction of the spiral-chain in
clockwise direction in $G_{ad}$ we see that for each edge-crossing the
direction of spiral-chain is switched from clockwise to anti-clockwise and 
\textit{vice versa. \ }Therefore spiral-chain edge-coloring algorithm
described in the previous section can be applied without any difficulty to
the graph $G_{ad}$.

For example illustration of apex and doublecross graphs are taken from [11]
and have been shown with the spiral-chain edge colorings.

\section{Conclusion}

\bigskip In this paper we have shown the effectiveness of spiral-chain
coloring in the proofs of the well-known Tait's and Tutte's three
edge-coloring of cubic graphs. With this result we have also given in a way
the proof of the four color theorem in two-ways either by using spiral-chain
vertex-coloring in the maximal planar or by using spiral-chain edge-coloring
for bridgless cubic planar graphs. We hope to announce similar results in an
near future for the Steinberg's vertex three-coloring conjecture of planar
graphs and for the Hadwiger's conjecture which is a generalization of the
four color problem.

\bigskip

\FRAME{ftbpFU}{248.75pt}{248.75pt}{0pt}{\Qcb{Spiral-chain coloring of an
uniquely edge 3-colorable cubic graph.}}{}{Figure}{\special{language
"Scientific Word";type "GRAPHIC";maintain-aspect-ratio TRUE;display
"USEDEF";valid_file "T";width 248.75pt;height 248.75pt;depth
0pt;original-width 531.875pt;original-height 531.875pt;cropleft "0";croptop
"1";cropright "1";cropbottom "0";tempfilename
'IIZPDG04.wmf';tempfile-properties "XPR";}}

\bigskip \FRAME{ftbpFU}{410.125pt}{196.5pt}{0pt}{\Qcb{Spiral-chain coloring
of the Tutte graph.}}{}{Figure}{\special{language "Scientific Word";type
"GRAPHIC";display "USEDEF";valid_file "T";width 410.125pt;height
196.5pt;depth 0pt;original-width 870.125pt;original-height 599pt;cropleft
"-0.3018";croptop "1";cropright "1.3018";cropbottom "0";tempfilename
'IIZPQA06.wmf';tempfile-properties "XPR";}}

\bigskip

\begin{center}
\bigskip \FRAME{ftbpFU}{223.125pt}{212.125pt}{0pt}{\Qcb{An apex and its
spiral-chain coloring.}}{}{Figure}{\special{language "Scientific Word";type
"GRAPHIC";maintain-aspect-ratio TRUE;display "USEDEF";valid_file "T";width
223.125pt;height 212.125pt;depth 0pt;original-width 644.75pt;original-height
612.6875pt;cropleft "0";croptop "1";cropright "1";cropbottom
"0";tempfilename 'IIZPB003.wmf';tempfile-properties "XPR";}}\FRAME{ftbpFU}{%
304.0625pt}{203pt}{0pt}{\Qcb{A doublecross and its spiral-chain coloring.}}{%
}{Figure}{\special{language "Scientific Word";type
"GRAPHIC";maintain-aspect-ratio TRUE;display "USEDEF";valid_file "T";width
304.0625pt;height 203pt;depth 0pt;original-width 637.3125pt;original-height
424.4375pt;cropleft "0";croptop "1";cropright "1";cropbottom
"0";tempfilename 'IIZP8L02.wmf';tempfile-properties "XPR";}}
\end{center}

\bigskip

\bigskip

\bigskip

\bigskip

\bigskip

\bigskip

\bigskip

\bigskip

\bigskip

{\Large References}

[1] \ \ P.G. Tait, "On the colouring of maps", \textit{Proceedings of the
Royal Society of Edinburgh Section A}, 10: 501-503, 1878-1880.

[2] \ \ P.G. Tait, "Remarks on the previous communication", \textit{%
Proceedings of the Royal Society of Edinburgh Section A}, 10:729, 1878-1880.

[3] \ \ P.G. Tait, "Note on atheorem in geometry of position", \textit{%
Transactions of the Royal Society of Edinburgh}, 29:657-660, 1880.

[4] \ \ A. Czumaj and A. Gibbons,"Guthrie's problem: new equivalances and
rapid reductions", \textit{Theoretical Computer Science}, 154:3-22, 1996.

[5] \ \ P. Sant, "Algorithmics of edge-colouring pairs of 3-regular trees",
Ph.D. thesis, University of London, December 2003.

[6]\ \ \ N. Robertson, P.D. Seymour and R.Thomas, "Tutte's edge-colouring
conjecture", \textit{J. Combin. Theory Ser. B}. 70 (1997), 166-183.

[7] \ \ N. Robertson, P.D. Seymour and R.Thomas, "Cyclically 5-connected
cubic graphs", in preparation.

[8] \ \ N. Robertson, P.D. Seymour and R.Thomas, "Excluded minors in cubic
graphs", in preparation.

[9] \ \ D. P. Sanders and R.Thomas, "Edge 3-coloring cubic apex graphs", in
preparation.

[10] D. P. Sanders, P.D. Seymour and R.Thomas, "Edge 3-coloring cubic
doublecross graphs", in preparation.

[11] W.T. Tutte,"On algebraic theory of graph coloring", \textit{J. Combin.
Theory}, 1 (1966), 15-50.

[12] R. Thomas, "An update on the four-color theorem", \textit{Notices Amer.
Math. Soc}.,45 (1998), 7, 848-859.

[13] K. Appel and W. Haken, "Every planar map is four colorable", \textit{%
Contemporary Math.}, 98 (1989).

[14] N. Robertson, D. P. Sanders, P. D. Seymour and R. Thomas,"The four
colour theorem",\textit{\ J. Combin. Theory Ser. B}. 70 (1997), 2-44.

[15] I. Cahit, "Spiral chains: A new proof of the four color theorem",arXiv
preprint, math.CO/0408247, August 18, 2004.

\bigskip 

\bigskip 

\bigskip 

I. Cahit

Near East University, North Cyprus

E-mail:\texttt{\ icahit@ebim.com.tr or icahit@gmail.com}

\bigskip 

July 6, 2005

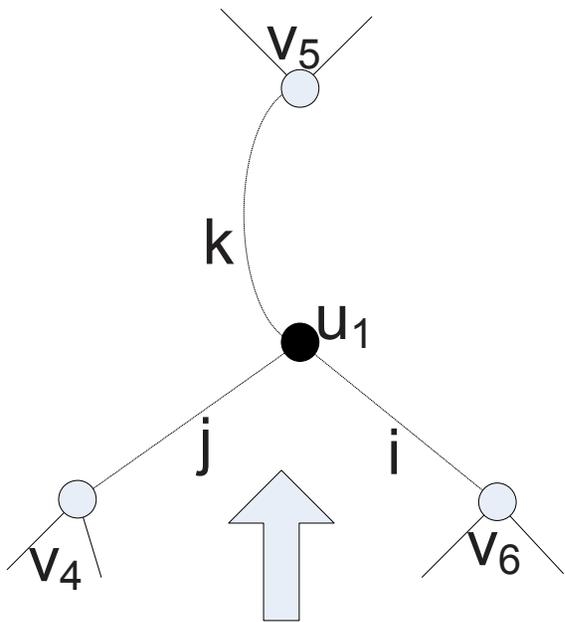
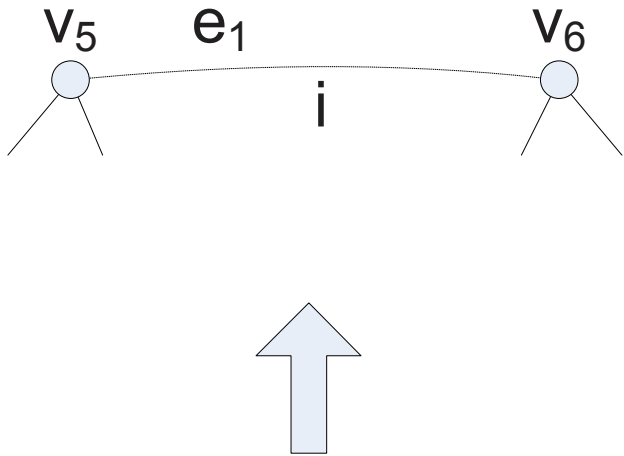
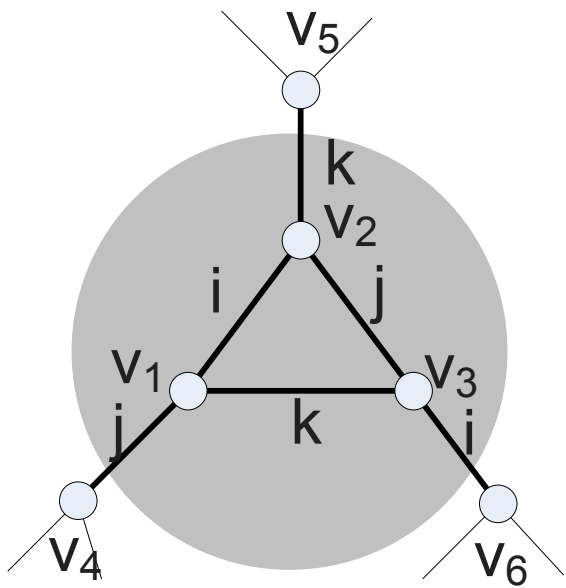
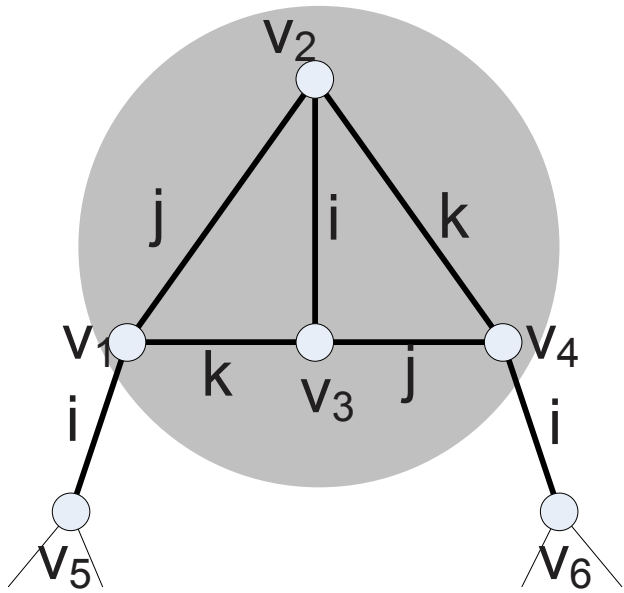
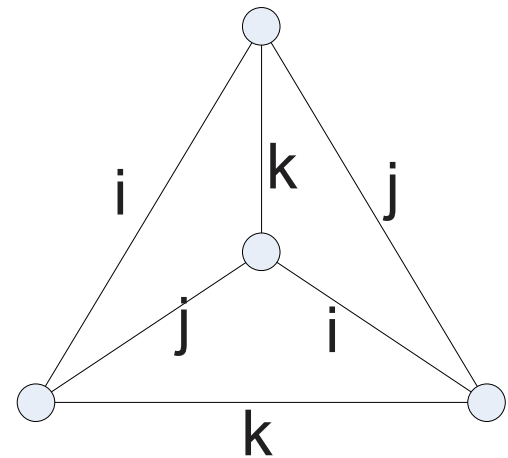
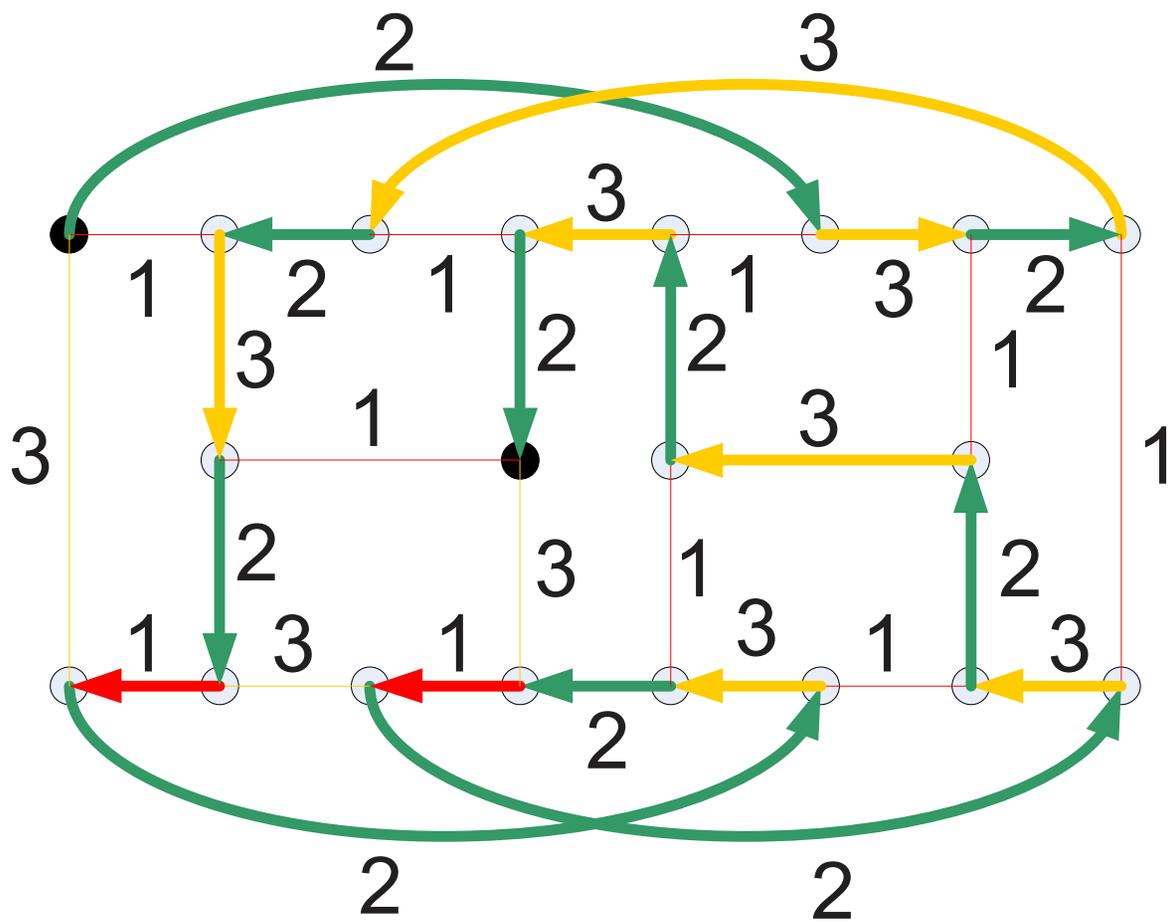
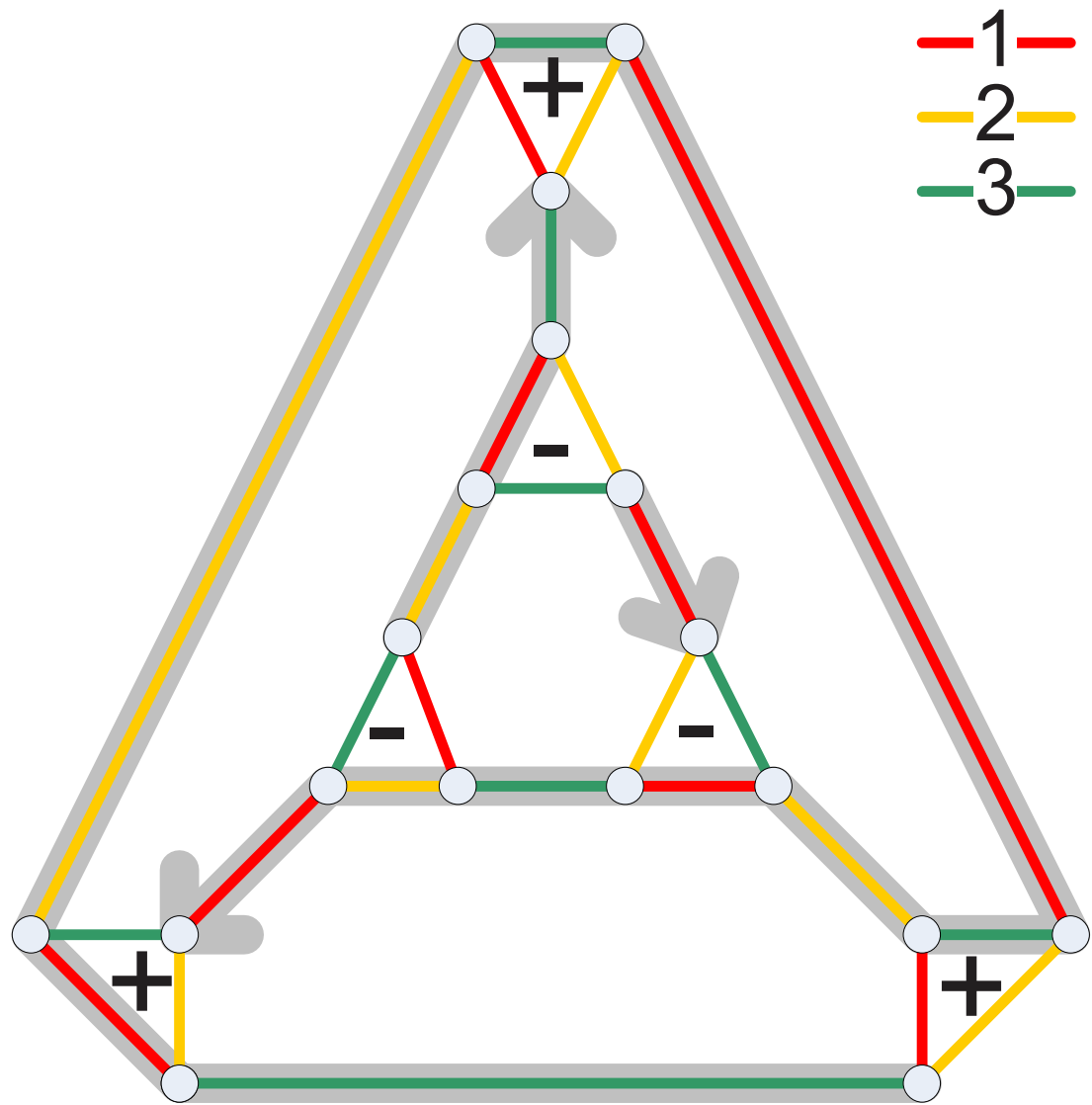
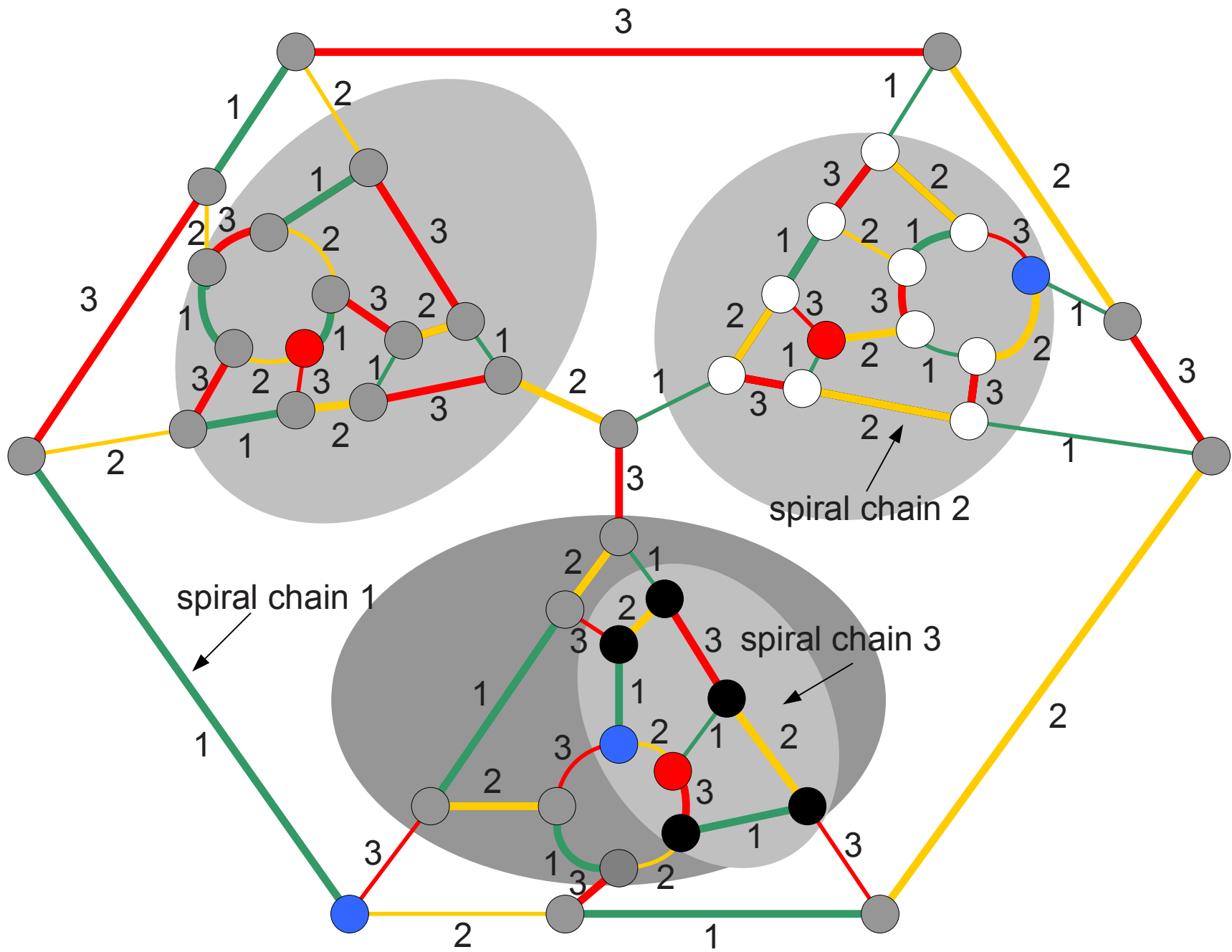
\end{document}